\documentclass{amsart}
\usepackage{amssymb}
\newtheorem{theorem}{Theorem}
\newtheorem{lemma}{Lemma}
\newtheorem{proposition}{Proposition}

\newtheorem{definition}{Definition}
\newtheorem{remark}{Remark}
\newcommand{\C}{\mathbb{C}}

\newcommand{\N}{\mathbb{N}}

\newcommand{\CP}{\mathbb{CP}}
\newcommand{\J}{\mathcal{J}}
\newcommand{\dbar}{\bar\partial}
\title[A remark about Donaldson's construction]{A remark about 
Donaldson's construction of symplectic submanifolds}
\author{D. Auroux}
\address{Department of Mathematics, M.I.T., Cambridge MA 02139, USA}
\email{auroux@math.mit.edu}

\begin{document}
\begin{abstract}
We describe a simplification of Donaldson's arguments for the construction 
of symplectic hypersurfaces \cite{D1} or Lefschetz pencils \cite{D2} that 
makes it possible to avoid any reference to Yomdin's work on the complexity
of real algebraic sets.
\end{abstract}

\maketitle

\section{Introduction}
Donaldson's construction of symplectic submanifolds \cite{D1} is
unquestionably one of the major results obtained in the past ten years
in symplectic topology. What sets it apart from many of the results obtained
during the same period is that it appeals neither to Seiberg-Witten theory,
nor to pseudo-holomorphic curves; in fact, most of Donaldson's argument is
a remarkable succession of elementary observations, combined in a
particularly clever way. One ingredient of the proof that does not qualify
as elementary, though, is an effective version of Sard's theorem for
approximately holomorphic complex-valued functions over a ball in $\C^n$
(Theorem 20 in \cite{D1}). The proof of this result, which occupies a
significant portion of Donaldson's paper (\S 4 and \S 5 of \cite{D1}),
appeals to very subtle considerations about the complexity of real algebraic
sets, following ideas of Yomdin \cite{Y}.

Methods similar to those in \cite{D1} were subsequently used to perform
various other constructions, leading in particular to Donaldson's
result that symplectic manifolds carry structures of symplectic 
Lefschetz pencils \cite{D2}, or to the result that symplectic 4-manifolds
can be realized as branched coverings of $\CP^2$ \cite{A2}. It was observed
in \cite{Aseoul} that, whereas Donaldson's construction of submanifolds 
can be thought of in terms of an estimated transversality result for sections
of line bundles, the subsequent constructions can be interpreted in terms of
estimated transversality with respect to stratifications in jet bundles.

As remarked at the end of \S 4 in \cite{Aseoul}, the transversality of the
$r$-jet of a section to a given submanifold in the bundle of $r$-jets is
equivalent to the non-intersection of the $(r+1)$-jet of the section with
a certain (possibly singular) submanifold of greater codimension in the
bundle of $(r+1)$-jets. This is of particular interest because the effective
Sard theorem for approximately holomorphic functions from $\C^n$ to $\C^m$
admits a conceptually much easier proof in the case where $m>n$ \cite{A2}.
In the case of the construction of symplectic submanifolds, the formalism
of jet bundles can be completely eliminated from the presentation; the
purpose of this note is to present the resulting simplified argument for
Donaldson's result (\S\S 2--3). We also observe (see \S 4) that a similar
simplification is possible for the higher-rank local result required for the
construction of symplectic Lefschetz pencils \cite{D2}.

\section{Overview of Donaldson's argument}

We first review Donaldson's construction of symplectic submanifolds
\cite{D1}, using the terminology and
notations of \cite{A2}. Let $(X^{2n},\omega)$ be a compact symplectic
manifold, and assume that the cohomology class $\frac{1}{2\pi}[\omega]$ is
integral. Endow $X$ with an $\omega$-compatible almost-complex structure $J$ 
and the corresponding Riemannian metric $g=\omega(\cdot,J\cdot)$. Consider a
Hermitian line bundle $L$ over $X$ such that $c_1(L)=\frac{1}{2\pi}[\omega]$,
equipped with a Hermitian connection $\nabla$ having curvature $-i\omega$.
The almost-complex structure $J$ induces a splitting of the connection into
$\nabla=\partial+\dbar$. We are interested in approximately holomorphic 
sections of the line bundles $L^{\otimes k}$ ($k\gg 0$)
satisfying a certain estimated transversality property: indeed, if we
can find a section $s$ such that $|\dbar s|\ll |\partial s|$ at every point 
where $s$ vanishes, then the zero set of $s$ is automatically a smooth 
symplectic submanifold in $X$ (cf.\ e.g.\ Proposition 3 of \cite{D1}). The
philosophical justification of the construction is that, as the twisting
parameter $k$ increases, one starts probing the geometry of $X$ at very
small scales, where the effects due to the non-integrability of $J$ become
negligible. This phenomenon is due to the curvature $-ik\omega$
of the connection on $L^{\otimes k}$, and leads us to work with a rescaled
metric $g_k=k\,g$ (the metric induced by $J$ and $k\omega$).

Let $(s_k)_{k\gg 0}$ be a sequence of sections of Hermitian vector bundles
$E_k$ equipped with Hermitian connections over $X$.
We make the following definitions:

\begin{definition}
The sections $s_k$ are {\em asymptotically holomorphic} if there exist 
constants $(C_p)_{p\in\N}$ such that, for all $k$ and at every
point of $X$, $|s_k|\le C_0$, $|\nabla^p s_k|_{g_k}\le C_p$ and 
$|\nabla^{p-1}\dbar s_k|_{g_k}\le C_p k^{-1/2}$ for all $p\ge 1$.
\end{definition}

\begin{definition}
The sections $s_k$ are {\em uniformly transverse to $0$} if there exists
a constant $\eta>0$ independent of $k$ such that the sections $s_k$ are
{\em $\eta$-transverse to $0$}, i.e.\ if at any point $x\in X$
where $|s_k(x)|<\eta$, the linear map $\nabla s_k(x):T_xX\to (E_k)_x$ 
is surjective and has a right inverse of norm less than $\eta^{-1}$ 
w.r.t.\ the metric $g_k$.
\end{definition}

When $\mathrm{rank}(E_k)>n$, uniform transversality means that
$|s_k(x)|\ge \eta$ at every point of $X$; on the other hand, when $E_k$ is
a line bundle and the sections $s_k$ are asymptotically holomorphic, uniform
transversality can be rephrased as a uniform lower bound on $|\partial s_k|$
at all points where $|s_k|<\eta$ (which by the above observation is enough
to ensure the symplecticity of $s_k^{-1}(0)$ for large $k$).
With this terminology, {\mbox Donaldson}'s result can be reformulated as 
follows (cf.\ Theorem 5 of \cite{D1}):

\begin{theorem}
For large values of $k$, the line bundles $L^{\otimes k}$ admit sections
$s_k$ that are asymptotically holomorphic and uniformly transverse to $0$.
\end{theorem}

The proof of Theorem 1 starts with a couple of preliminary lemmas about
the existence of approximately holomorphic rescaled Darboux
coordinates on $X$ and of large families of well-concentrated asymptotically 
holomorphic sections of $L^{\otimes k}$.

\begin{lemma}
There exists a constant $c>0$ such that
near any point $x\in X$, for any integer $k$, there exist local complex 
Darboux coordinates $z_k=(z_k^1,\dots,z_k^n):(X,x)\to(\C^n,0)$ for the 
symplectic structure $k\omega$, such that the following estimates hold 
uniformly in $x$ and $k$ at every point of the ball $B_{g_k}(x,c\sqrt{k})$:
$|z_k(y)|=O(\mathrm{dist}_{g_k}(x,y))$, $|\dbar z_k(y)|_{g_k}=O(k^{-1/2}
\mathrm{dist}_{g_k}(x,y))$, $|\nabla^r\dbar z_k|_{g_k}=O(k^{-1/2})$,
$|\nabla^r z_k|_{g_k}=O(1)$ $\forall r\ge 1$;
and denoting by $\psi_k:(\C^n,0)\to (X,x)$ the inverse map, the estimates
$|\dbar\psi_k(z)|_{g_k}=O(k^{-1/2}|z|)$,
$|\nabla^r\dbar\psi_k|_{g_k}=O(k^{-1/2})$
and $|\nabla^r\psi_k|_{g_k}=O(1)$ hold $\forall r\ge 1$
at every point of the ball $B_{\C^n}(0,c\sqrt{k})$,
where $\dbar\psi_k$ is defined with respect to the almost-complex structure 
$J$ on $X$ and the standard complex structure on $\C^2$.
\end{lemma}

Lemma 1 is identical to Lemma 3 of \cite{A2}, or to the discussion on pp.\
674--675 of \cite{D1} if one keeps track carefully of the available
estimates; the idea is simply to start with usual Darboux coordinates for
$\omega$, compose them with a linear transformation to ensure holomorphicity
at the origin, and then rescale them by a factor of $\sqrt{k}$.

\begin{definition}
A section $s$ of $E_k$ has {\em Gaussian decay} in $C^r$ 
norm away from a point $x\in X$ if there exist a polynomial $P$ and a 
constant $\lambda>0$ such that for all $y\in X$, $|s(y)|$, 
$|\nabla s(y)|_{g_k}$, \ldots, $|\nabla^r s(y)|_{g_k}$ are all bounded 
by $P(d(x,y))\exp(-\lambda\,d(x,y)^2)$, where $d(.,.)$ is the distance
induced by $g_k$.
The decay properties of a family of sections are said to be {\em uniform}
if $P$ and $\lambda$ can be chosen independently of $k$ and of the point 
$x$ at which decay occurs for a given section.
\end{definition}

\begin{lemma}
Given any point $x\in X$, for all large enough $k$, there exist
asymptotically holomorphic sections  $s_{k,x}^\mathrm{ref}$ of 
$L^{\otimes k}$ over $X$, such that
$|s_{k,x}^\mathrm{ref}|\ge c_0$ at every point of the ball of 
$g_k$-radius $1$ centered at $x$, for some universal constant $c_0>0$, and
such that the sections $s_{k,x}^\mathrm{ref}$
have uniform Gaussian decay away from $x$ in $C^2$ norm.
\end{lemma}

Lemma 2 is essentially Proposition 11 of \cite{D1}. Considering
a local trivialization of $L^{\otimes k}$ where the connection 1-form 
is $\frac{1}{4}\sum (z_k^j d\bar{z}_k^j - \bar{z}_k^j dz_k^j)$, the
sections $s_{k,x}^\mathrm{ref}$ are constructed by multiplication of
the function
$\exp(-|z_k|^2/4)$ by a suitable cut-off function at
distance $k^{1/6}$ from the origin.

The central ingredient is the following result about the near-critical sets
of approximately holomorphic functions (used in the special case $m=1$):

\begin{proposition}
Let $f$ be a function defined over the ball $B^+$ of radius $\frac{11}{10}$
in $\C^n$ with values in $\C^m$. Let $\delta$ be a constant with
$0<\delta<\frac{1}{4}$, and let $\eta=\delta\log(\delta^{-1})^{-p}$ where
$p$ is a fixed integer depending only on $n$ and $m$. Assume
that $f$ satisfies the bounds $|f|_{C^0(B^+)}\le 1$ and $|\dbar
f|_{C^1(B^+)}\le \eta$. Then there exists $w\in\C^m$ with $|w|\le \delta$
such that $f-w$ is $\eta$-transverse to $0$ over the interior ball $B$ of
radius~$1$.
\end{proposition}

The case $m=1$ is Theorem 20 of \cite{D1}; the comparatively much easier
case $m>n$ is Proposition 2 of \cite{A2}; the general case is proved in
\cite{D2}. In all cases the proof begins with an approximation of $f$ first
by a holomorphic function (using general elliptic theory), then by a
polynomial $g$ of degree $O(\log(\eta^{-1}))$ (by truncating the power series
expansion at the origin). The proof in the case $m=1$ then appeals to a
rather sophisticated result on the complexity of real algebraic sets to 
control the size of the set of points where $\partial g$ is small (the
near-critical points) \cite{D1}. Meanwhile, in the case $m>n$, since we
only have to find $w$ such that $|f-w|\ge \eta$ at every point of $B$,
it is sufficient to observe that the image of the polynomial map $g$
is contained in a complex algebraic hypersurface $H$ in $\C^m$; the result then
follows from a standard result about the volume of a tubular neighborhood
of $H$, which can be estimated using an explicit bound on the degree of $H$
\cite{A2}.

Given asymptotically holomorphic sections $s_k$ of $L^{\otimes k}$ and a
point $x\in X$, one can apply Proposition 1 to the complex-valued functions 
$f_k=s_k/s_{k,x}^\mathrm{ref}$ (defined over
a neighborhood of $x$) in order to find constants $w_k$ such that
the functions $f_k-w_k$ are uniformly transverse to $0$ near $x$;
multiplying by $s_{k,x}^\mathrm{ref}$, it follows that the sections
$s_k-w_ks_{k,x}^\mathrm{ref}$ are uniformly transverse to $0$ near $x$.
Therefore, we have:

\begin{proposition}
There exist constants $c$, $c'$, $p$, $\delta_0>0$ such that, given a
real number $\delta\in (0,\delta_0)$, a sequence of asymptotically 
holomorphic sections $s_k$ of $L^{\otimes k}$ and a point $x\in X$,
for large enough $k$ there exist asymptotically holomorphic sections
$\tau_{k,x}$ of $L^{\otimes k}$ with the following properties: {\rm (a)}
$|\tau_{k,x}|_{C^1,g_k}<\delta$, {\rm (b)} the sections 
$\frac{1}{\delta}\tau_{k,x}$ have uniform Gaussian
decay away from $x$ in $C^1$ norm, and {\rm (c)} the sections 
$s_k+\tau_{k,x}$ are $\eta$-transverse to $0$ at every point of the ball 
$B_{g_k}(x,c)$, with $\eta=c'\delta\log(\delta^{-1})^{-p}$.
\end{proposition}

This result lets us achieve estimated transversality over a small ball in $X$
by adding to $s_k$ a small well-concentrated perturbation. Uniform 
transversality over the entire manifold $X$ is achieved by proceeding 
iteratively, adding successive perturbations to the sections in order to 
obtain transversality properties over larger and larger subsets of $X$. 
The key observation is that estimated transversality is an open property 
(preserved under $C^1$-small perturbations). Since the transversality estimate 
decreases after each perturbation, it is important to obtain global
uniform transversality after a number of steps that remains bounded
independently of $k$; this is made possible by the uniform decay properties
of the perturbations, using a beautiful observation of Donaldson.
The reader is referred to \S 3 of \cite{D1} or to Proposition 3 of \cite{A2}
for details.

\section{The simplified argument}

Keeping the same general strategy, the proof of Theorem 1 can be simplified
by appealing to a result weaker than Proposition 1, namely the following
statement:

\begin{proposition}
Let $f$ be a function defined over the ball $B^+$ of radius $\frac{11}{10}$
in $\C^n$ with values in $\C$. Let $\delta$ be a constant with
$0<\delta<\frac{1}{4}$, and let $\eta=\delta\log(\delta^{-1})^{-p'}$ where
$p'$ is a fixed integer depending only on $n$. Assume
that $f$ satisfies the bounds $|f|_{C^1(B^+)}\le 1$ and $|\dbar
f|_{C^2(B^+)}\le \eta$. Then there exists $w=(w_0,w_1,\dots,w_n)\in\C^{n+1}$ 
with $|w|\le \delta$ such that the function $f-w_0-\sum w_i z_i$ is 
$\eta$-transverse to $0$ over the interior ball $B$ of radius~$1$.
\end{proposition}

\proof
Let $g=(g_0,\dots,g_n):B^+\to \C^{n+1}$ be the function defined by
$g_i=\partial f/\partial z_i$ for $1\le i\le n$ and $g_0=f-\sum_{i=1}^n 
z_i g_i$. The bounds on $f$ immediately imply that $|g|_{C^0(B^+)}\le C_n$
and $|\dbar g|_{C^1(B^+)}\le C_n\eta$, for some constant $C_n$ depending
only on the dimension. We can safely choose the constant $p'$ appearing in
the definition of $\eta$ to be larger than the constant $p$ appearing in
Proposition 1. Therefore we can apply Proposition 1 in 
its easy version ($m=n+1$) to the function $g$, after scaling down by the 
constant factor $C_n$. This gives us a constant $w=(w_0,\dots,w_n)
\in\C^{n+1}$, bounded by $\delta$, and such that $|g-w|\ge \alpha$ at 
every point of $B$, where $\alpha=\delta\log((\delta/C_n)^{-1})^{-p}$.

Define $\tilde{f}=f-w_0-\sum w_i z_i$ and $\tilde{g}=g-w$, and observe that
$\partial \tilde{f}/\partial z_i=\tilde{g}_i$ for $1\le i\le n$ and
$\tilde{f}-\sum_{i=1}^n z_i\tilde{g}_i=\tilde{g}_0$. Let $z\in B$ be
a point where $|\partial \tilde{f}|<\frac{1}{4}\alpha$. Since 
$\partial \tilde{f}/\partial z_i=\tilde{g}_i$, and since $|\tilde{g}(z)|\ge
\alpha$ by construction, we have the inequality $|\tilde{g}_0(z)|>
\frac{3}{4}\alpha$. However, $|\tilde{f}(z)-\tilde{g}_0(z)|=|\sum z_i
\tilde{g}_i(z)|\le |z|\,|\partial \tilde{f}(z)|<\frac{1}{4}\alpha$
(recall that $z$ belongs to the unit ball). Therefore
$|\tilde{f}(z)|>\frac{1}{2}\alpha$. 

Conversely, at any point $z\in B$ where
$|\tilde{f}|\le \frac{1}{2}\alpha$ we must have $|\dbar\tilde{f}(z)|\ge
\frac{1}{4}\alpha$.
However, because of the bound on $\dbar \tilde{f}=\dbar f$, if we
assume that $\eta<\frac{1}{8}\alpha$ then this inequality implies that
$\nabla \tilde{f}(z)$ is surjective and admits a right inverse of norm 
at most $(\frac{1}{8}\alpha)^{-1}$. Hence we conclude from the previous 
discussion that $\tilde{f}$ is $\frac{1}{8}\alpha$-transverse to $0$ over
$B$. Finally, we observe that, because $\delta<\frac{1}{4}$, if the constant
$p'$ is chosen large enough then $\eta=\delta\log(\delta^{-1})^{-p'}<
\frac{1}{8}\alpha=\frac{1}{8}\delta\log((\delta/C_n)^{-1})^{-p}$, so that
$\tilde{f}$ is $\eta$-transverse to $0$ over $B$.
\endproof

Although it is weaker, Proposition 3 is in practice
interchangeable with the case $m=1$ of Proposition 1, in particular for
the purpose of proving Proposition 2.

\proof[Proof of Proposition 2]
We use the same argument as Donaldson \cite{D1}: we work in approximately
holomorphic Darboux coordinates on a neighborhood of the given point $x$,
using Lemma 1. Using the sections $s_{k,x}^\mathrm{ref}$ given 
by Lemma 2 to define local trivializations of $L^{\otimes k}$, the
sections $s_k$ can be identified with complex-valued functions
$f_k=s_k/s_{k,x}^\mathrm{ref}$. The estimates on $s_k$ and
$s_{k,x}^\mathrm{ref}$ imply that the functions $f_k$ are approximately
holomorphic near the origin (in particular $|\dbar f_k|_{C^2}=O(k^{-1/2})$);
after a suitable rescaling of the coordinates and of the functions 
by uniform constant factors, we can assume additionally that 
$|f_k|_{C^1}\le 1$ near the origin, and that the estimates hold over a 
neighborhood of
the origin that contains the ball $B^+$. Therefore, the assumptions of 
Proposition 3 are satisfied provided that $k$ is sufficiently large
to ensure that $k^{-1/2}\ll \eta$.

By Proposition 3, we can find $w_k=(w_{k,0},\dots,w_{k,n})\in\C^{n+1}$, with
$|w_k|\le \delta$, such that $\tilde{f_k}=f_k-w_{k,0}-\sum w_{k,i}z_i$ is
$\gamma$-transverse to $0$ over the unit ball, where
$\gamma=\delta\log(\delta^{-1})^{-p'}$. Define $\tau_{k,x}=-w_{k,0}
s_{k,x}^\mathrm{ref}-\sum w_{k,i}z_k^is_{k,x}^\mathrm{ref}$. The estimates
on $z_k^i$ from Lemma 1 and on $s_{k,x}^\mathrm{ref}$ from Lemma 2 imply
that the sections $z_{k,i}s_{k,x}^\mathrm{ref}$ of $L^{\otimes k}$ are
asymptotically holomorphic and have uniform Gaussian decay away from $x$.
Therefore, it is easy to check that the sections $\frac{1}{\delta}\tau_{k,x}$
are asymptotically holomorphic and have uniform Gaussian decay. Moreover,
because there exist uniform bounds on $s_{k,x}^\mathrm{ref}$ and $z_k^i
s_{k,x}^\mathrm{ref}$, one easily checks that $|\tau_{k,x}|_{C^1,g_k}$ is
bounded by some constant multiple of $\delta$; decreasing the required 
bound on $|w_k|$, we can assume that the constant is equal to $1$, to the 
expense of inserting a constant factor in the above expression for $\gamma$.
Finally, observing that $s_k+\tau_{k,x}=\tilde{f}_k s_{k,x}^\mathrm{ref}$
over a neighborhood of $x$, it is straightforward to check that the
$\gamma$-transversality to $0$ of $\tilde{f}_k$ and the lower bound
satisfied by $s_{k,x}^\mathrm{ref}$ imply a uniform transversality 
property of the desired form for the section $s_k+\tau_{k,x}$.
\endproof

\begin{remark} \rm
Proposition 3 also admits a version for one-parameter families of functions:
given functions $f_t:B^+\to \C$ depending continuously on a parameter $t\in
[0,1]$ and satisfying the assumptions of Proposition 3 for all values of
$t$, we can find constants $w_t\in\C^{n+1}$, depending continuously on $t$,
such that the conclusion holds for all values of $t$. This is because the
auxiliary functions $g_t:B^+\to\C^{n+1}$ introduced in the proof also
depend continuously on $t$, which allows us to appeal to the one-parameter 
version of Proposition 1 (cf.\ e.g.\ Proposition 2 of \cite{A2}). We can
therefore simplify the argument proving the asymptotic uniqueness of
the constructed submanifolds \cite{A1} in the same manner as the
construction itself.
\end{remark}

\begin{remark} \rm
The idea behind the modified argument can be interpreted as follows in terms
of $1$-jets of sections: let $\J^1 L^{\otimes k}=L^{\otimes k}\oplus
(T^*X^{1,0}\otimes L^{\otimes k})$, and define the $1$-jet of a section
$s_k\in\Gamma(L^{\otimes k})$ as $j^1s_k=(s_k,\partial
s_k)\in\Gamma(\J^1 L^{\otimes k}$). The jet bundles carry natural
Hermitian metrics (induced by those on $L^{\otimes k}$ and the metrics $g_k$
on the cotangent bundle), and natural Hermitian connections for which
the $1$-jets of asymptotically holomorphic sections of $L^{\otimes k}$ are
asymptotically holomorphic sections of $\J^1 L^{\otimes k}$. It is worth
noting that the natural connection on $\J^1 L^{\otimes k}$ is {\em not} 
the connection $\nabla$ induced by the connection on $L^{\otimes k}$ and 
the Levi-Civita connection, because $\dbar^\nabla(s_k,\partial s_k)=
(\dbar s_k,\dbar \partial s_k)$ differs
from $(\dbar s_k,-\partial\dbar s_k)$ (which is bounded by $O(k^{-1/2})$) 
by the curvature term $-ik\omega s_k$. Therefore, we must instead work with
the Hermitian connection $\tilde\nabla$ characterized by the
formula $\dbar^{\tilde\nabla}(\sigma^0,\sigma^1)=
\dbar^\nabla (\sigma^0,\sigma^1)+(0,ik\omega\sigma^0)$, where $\omega$ is
viewed as a $(0,1)$-form with values in $T^*X^{1,0}$.

Observe that the 1-jets $j^1\sigma_{k,x,0},\dots,j^1\sigma_{k,x,n}$, where
$\sigma_{k,x,0}=s_{k,x}^\mathrm{ref}$ and $\sigma_{k,x,i}=z_k^i
s_{k,x}^\mathrm{ref}$ for $1\le i\le n$, are asymptotically holomorphic
sections of $\J^1 L^{\otimes k}$, with uniform Gaussian decay away from $x$,
which form a local frame of the jet bundle over a neighborhood of $x$. 
Therefore, given asymptotically holomorphic sections $s_k$ and a point 
$x\in X$, there exist local complex-valued functions $g_{k,0},\dots,g_{k,n}$
such that $j^1 s_k=\sum g_{k,i}\, j^1 \sigma_{k,x,i}$. Moreover, remark
that a section of $L^{\otimes k}$ is uniformly transverse to $0$ if and only
if its $1$-jet satisfies a uniform lower bound. Therefore, our argument
actually amounts to a local perturbation of $j^1 s_k$, using the given
local frame $\{j^1\sigma_{k,x,i}\}$, in order to bound it away from $0$;
because the rank of the jet bundle is $n+1>n$, the easy version of
Proposition 1 is sufficient for that purpose.
The curious reader is referred to \cite{Aseoul} for a more detailed
discussion of estimated transversality using the formalism of jet bundles.
\end{remark}

\section{The higher-rank local result}

We now formulate and prove an analogue of Proposition 3 for functions
with values in $\C^m$ ($m\le n$); as in the case $m=1$, the statement 
differs from Proposition 1 by allowing the extra freedom of affine 
perturbations rather than restricting oneself to constants.

\begin{proposition}
Let $f$ be a function defined over the ball $B^+$ of radius $\frac{11}{10}$
in $\C^n$ with values in $\C^m$, $m\le n$. Let $\delta$ be a constant with
$0<\delta<\frac{1}{4}$, and let $\eta=\delta\log(\delta^{-1})^{-p'}$ where
$p'$ is a fixed integer depending only on $m$ and $n$. Assume
that $f$ satisfies the bounds $|f|_{C^0(B^+)}\le 1$ and $|\dbar
f|_{C^1(B^+)}\le \eta$. Then there exists 
$w=(w_0,w_1,\dots,w_n)\in\C^{m(n+1)}$ $($each $w_i$ is an element of $\C^m)$
with $|w|\le \delta$ such that the function $f-w_0-\sum w_i z_i$ is 
$\eta$-transverse to $0$ over the interior ball $B$ of radius~$1$.

Moreover, given a one-parameter family of functions $f_t:B^+\to \C$ depending 
continuously on a parameter $t\in [0,1]$ and satisfying the above assumptions
for all $t$, we can find constants $w_t\in\C^{m(n+1)}$, depending continuously
on $t$, such that the conclusion holds for all values of $t$.
\end{proposition}

This statement is essentially interchangeable with Proposition 1 for all
practical applications, and in particular the case $m=n$ allows us to
simplify noticeably the argument for Donaldson's construction of symplectic
Lefschetz pencils \cite{D2}. Indeed, the main problem to be solved
is the following: given pairs of asymptotically holomorphic sections
$(s_k^0,s_k^1)$ of $L^{\otimes k}$, defining $\CP^1$-valued maps 
$f_k=[s_k^0:s_k^1]$ away from the base loci, one must perturb them so that
the differentials $\partial f_k$ (which are sections of rank $n$ vector 
bundles) become uniformly transverse to $0$. This ensures the non-degeneracy 
of the singular points of the pencil. The manner in which the problem
reduces to the $m=n$ case of Proposition 1 is explained in detail in
\cite{D2}, and the reduction to Proposition 4 is essentially identical
except that the resulting perturbations of $(s_k^0,s_k^1)$ are products
of $s_{k,x}^\mathrm{ref}$ by quadratic (rather than linear) polynomials.

\proof
Although for technical reasons we cannot use directly the case $m>n$ of
Proposition 1, the argument presents many similarities with \S 2.3 of
\cite{A2}; we accordingly skip the details whenever the two arguments
parallel each other in an obvious manner.
As in the case of Proposition 1, we first use the bounds on $f$ to find an
approximation by a polynomial $h:\C^n\to \C^m$ of degree 
$d=O(\log(\eta^{-1}))$ such that $|h-f|_{C^1(B)}\le c\,\eta$ for some
constant $c$ (see Lemmas 27 and 28 of \cite{D1}). Observe that, if we can
perturb $h$ by less than $\delta$ to make it $(c+1)\eta$-transverse to $0$
over $B$, then because transversality is an open property the desired result
on $f$ will follow immediately. So we are reduced to the case of
a polynomial function $h=(h^1,\dots,h^m)$ of degree $d=O(\log(\eta^{-1}))$.

If $w=(w_0,\dots,w_n)$ is a vector in $\C^{m(n+1)}$, denote by 
$(w_i^j)_{1\le j\le m}$ the components of $w_i$, and let
$\vec{w}=(w_1,\dots,w_n)\in \C^{m\times n}$.
The set of choices to be avoided 
for $w$ is $$S=\{w\in \C^{m(n+1)},\ \exists z\in B\ \mbox{s.t.}
\ h(z)-w_0-{\textstyle\sum} w_i z_i=0,\ {\textstyle \bigwedge^m}(\partial 
h(z)-\vec{w})=0\}.$$ Indeed, observe that $h-w_0-\sum w_i z_i$ is transverse
to $0$ over $B$ (without any estimate) if and only if $w\not\in S$. We now
define a polynomial function $g:\C^{N-1}\to\C^N$, where $N=m(n+1)$, which
parametrizes a dense subset of $S$. Given elements $z=(z_i)_{1\le i\le n}
\in\C^n$, $\theta=(\theta_i^j)_{1\le i\le n,\ 1\le j\le m-1}\in \C^{(m-1)n}$ and
$\lambda=(\lambda_j)_{1\le j\le m-1}\in \C^{m-1}$, we define
$g(z,\theta,\lambda)\in \C^{m(n+1)}$ by the formulas
$$\begin{cases}
g_i^j(z,\theta,\lambda)=\dfrac{\partial h^j}{\partial z_i}(z)+\theta_i^j & 
\mathrm{for}\ 1\le i\le n,\ 1\le j\le m-1,\\
g_i^m(z,\theta,\lambda)={\dfrac{\partial h^m}{\partial z_i}(z)+
\sum\limits_{j=1}^{m-1} \lambda_j\theta_i^j} & \mathrm{for}\ 1\le i\le n,\\
g_0^j(z,\theta,\lambda)=h^j(z)-\sum\limits_{i=1}^n g_i^j(z,\theta,\lambda) 
z_i & \mathrm{for}\ 1\le j\le m.
\end{cases}$$
One easily checks that the image by $g$ of the subset 
$\{(z,\theta,\lambda)\in\C^{N-1},\ z\in B\}$ is contained in $S$, in which
it is a dense subset. Observe that $g$ is a polynomial map with the same
degree $d$ as $h$ (provided that $d\ge 2$). Therefore, the image
$g(\C^{N-1})$ is contained in an algebraic surface $H\subset \C^N$, of
degree at most $D=N\,d^{N-1}$. Indeed, denoting by $E$ the space of
polynomials of degree at most $D$ in $N$ variables and by $E'$ the space of
polynomials of degree at most $dD$ in $N-1$ variables, we have $\dim E=
\binom{D+N}{N}>\binom{dD+N-1}{N-1}=\dim E'$, so that the map from $E$ to
$E'$ defined by $P\mapsto P\circ g$ cannot be injective, and a non-zero
element of its kernel provides an equation for the hypersurface $H$
(see \S 2.3 of \cite{A2} for details). 

Since $g(B\times\C^{(m-1)n}\times \C^{m-1})$ is dense in $S$, we conclude
that $S\subset H$. From this point on, the argument is very similar to 
\S 2.3 of \cite{A2}, to which the reader is referred for details. Standard
results on complex algebraic hypersurfaces, essentially amounting to the
well-known monotonicity formula, allow us to bound the size of $S$ and of 
its tubular neighborhoods (cf.\ e.g.\ Lemma 4 of \cite{A2}). 
In particular, denoting by $\bar{B}$ the ball of radius $\delta$ centered at
the origin in $\C^N$ and by $V_0$ the volume of the unit ball in dimension
$2N-2$, we have $vol_{2N-2}(H\cap \bar{B})\le DV_0 \delta^{2N-2}$, while
given any point $x\in H$ we have $vol_{2N-2}(H\cap B(x,\eta))\ge V_0
\eta^{2N-2}$. Therefore, choosing a suitable covering of $\bar{B}$ by balls
of radius $\eta$, one can show that $H\cap \bar{B}$ is contained in the
union of $M=C\,D\,\delta^{2N-2}\eta^{-(2N-2)}$ balls of radius $\eta$, where
$C$ is a constant depending only on $N$. As a consequence, the neighborhood
$Z=\{w\in\C^N,\ |w|\le\delta,\ \mathrm{dist}(w,S)\le (3c+3)\eta\}$ is
contained in the union of $M$ balls of radius $(3c+4)\eta$.

A simple comparison of the volumes implies that, if the constant $p'$ is
chosen suitably large, then the volume of $Z$ is much smaller than that of
the ball $\bar{B}$, and therefore $\bar{B}-Z$ is not empty, i.e.\ $\bar{B}$
contains an element $w$ which lies at distance more than $(3c+3)\eta$ from
$S$. Moreover, using a standard isoperimetric inequality we can show that 
$\bar{B}-Z$ contains a unique large connected component; it follows that,
in the case where the data depends continuously on a parameter $t\in [0,1]$,
the subset $\bigsqcup \{t\}\times (\bar{B}-Z_t)\subset [0,1]\times \bar{B}$ 
contains a preferred large connected component, in which we can choose
elements $w_t$ depending continuously on $t$.

To complete the proof of Proposition 4, we only need to show that, 
if $w\in \bar{B}$ lies at
distance more than $(3c+3)\eta$ from $S$, then $\tilde{f}=f-w_0-\sum w_iz_i$
is $\eta$-transverse to $0$ over $B$. In fact, it is sufficient to show
that $\tilde{h}=h-w_0-\sum w_iz_i$ is $(c+1)\eta$-transverse to $0$ over
$B$, because $|\tilde{h}-\tilde{f}|_{C^1(B)}=|h-f|_{C^1(B)}\le c\eta$ and
transversality is an open property. We conclude using the following lemma:

\begin{lemma}
If $w$ lies at distance more than $3\alpha$ from $S$ for some constant
$\alpha>0$, then $\tilde{h}=h-w_0-\sum w_iz_i$ is $\alpha$-transverse to 
$0$ over $B$.
\end{lemma}

To prove Lemma 3, we first provide an alternative definition of
$\alpha$-transversality:

\begin{lemma}
Let $L:E\to F$ be a linear map between Hermitian complex vector spaces, 
and choose a constant $\alpha>0$. The two following properties are 
equivalent:

$(i)$ $L$ is surjective and has a right inverse $R:F\to E$ of norm at most
$\alpha^{-1}$,

$(ii)$ for every unit vector $v$ in $F$, the component $\langle v,L\rangle
=v^*L$ of $L$ along $v$ is a linear form on $E$ such that $|v^*L|\ge\alpha$.
\end{lemma}

\proof
If $(i)$ holds, then given any unit vector $v\in F$, the vector $u=Rv$
is such that $|u|\le\alpha^{-1}$ and $\langle v,Lu\rangle=|v|^2=1$. Therefore
the linear form $\langle v,L\rangle$ has norm at least $\alpha$, and
$(ii)$ holds.

Conversely, assume $(ii)$ holds. Then for any $v\in F$ we have $|v^*L|\ge
\alpha |v|$, i.e.\ $v^*LL^*v\ge \alpha^2 |v|^2$. Therefore, the Hermitian
endomorphism $LL^*$ of $F$ is positive definite and has eigenvalues $\ge
\alpha^2$. It follows that it admits an inverse $U=(LL^*)^{-1}$ of operator
norm at most $\alpha^{-2}$. We have $LL^*U=\mathrm{Id}$, and $|L^*Uv|^2=
\langle v,ULL^*Uv\rangle=\langle v,Uv\rangle\le \alpha^{-2}|v|^2$, so that
$R=L^*U$ is a right inverse of norm at most $\alpha^{-1}$.
\endproof

\proof[Proof of Lemma 3]
Assume that $\tilde{h}$ is not $\alpha$-transverse to $0$ over $B$: using the
definition and Lemma 4, there exists a point $z\in B$ and a unit vector
$v\in \C^m$ such that $|\tilde{h}(z)|<\alpha$ and $|\langle v,\partial
\tilde{h}(z)\rangle|<\alpha$. Let $u=(u_0,u_1,\dots,u_n)\in 
\C^{m(n+1)}$ be such that $u_i=\langle v,
\partial\tilde{h}/\partial z_i\rangle\,v$ and $u_0=\tilde{h}(z)-\sum
z_i u_i$. We clearly have $|(u_1,\dots,u_n)|<\alpha$, and $|u_0|<2\alpha$,
so that $|u|<3\alpha$. On the other hand, if we consider the function
$\hat{h}=h-(w_0+u_0)-\sum (w_i+u_i)z_i$, then by construction $\hat{h}(z)=0$
and $\langle v,\partial \hat{h}(z)\rangle=0$. Therefore $w+u\in S$, 
and so $w$ is within distance $3\alpha$ of $S$.
\endproof

\end{document}